\documentclass[11pt]{article}
\usepackage{xspace,amssymb,amsmath,helvet,graphicx}
\begin{document}
\newcommand{\dee}{\,\mbox{d}}
\newcommand{\naive}{na\"{\i}ve }
\newcommand{\eg}{e.g.\xspace}
\newcommand{\ie}{i.e.\xspace}
\newcommand{\pdf}{pdf.\xspace}
\newcommand{\etc}{etc.\@\xspace}
\newcommand{\PhD}{Ph.D.\xspace}
\newcommand{\MSc}{M.Sc.\xspace}
\newcommand{\BA}{B.A.\xspace}
\newcommand{\MA}{M.A.\xspace}
\newcommand{\role}{r\^{o}le}
\newcommand{\signoff}{\hspace*{\fill} Rose Baker \today}
\newenvironment{entry}[1]%
{\begin{list}{}{\renewcommand{\makelabel}[1]{\textsf{##1:}\hfil}%
\settowidth{\labelwidth}{\textsf{#1:}}%
\setlength{\leftmargin}{\labelwidth}
\addtolength{\leftmargin}{\labelsep}
\setlength{\itemindent}{0pt}
}}%
{\end{list}}
\title{Creating new distributions using integration and summation by parts}

\author{Rose Baker, School of Business\\University of Salford, UK\\email r.d.baker@salford.ac.uk}
\maketitle
\begin{abstract}
Methods for generating new distributions from old can be thought of as techniques for simplifying integrals used in reverse.
Hence integrating a probability density function (pdf) by parts provides a new way of modifying distributions; the resulting pdfs are integrals that sometimes
require computation as special functions. 
Summation by parts can be used similarly for discrete distributions.
The general methodology is given, with some examples of distribution classes and of specific distributions, and  fits to data.
\end{abstract}

\section*{Keywords}
Mixture distribution; partial integration; stochastic dominance; summation by parts; discrete distribution; special functions
\section{Introduction}
Parametric models of probability distributions are essential for statistical inference. Hence a vast number of distributions of all types has been created,
and a common way to generate new distributions is to modify an old one. Jones (2015) reviews the main techniques for generalizing univariate symmetric distributions,
and Lai (2012) gives a comprehensive account of ways to modify survival distributions.

Transforming the random variable is probably the most popular method. 
It is a technique that is often used to evaluate unknown integrals, and is used `in reverse' where the integral of the pdf (unity) is already known, to generate new pdfs.
For example, the exponential distribution with survival function $\bar{F}(x)=\exp(-\alpha x)$
becomes the Weibull distribution with $\bar{F}(y)=\exp(-(\alpha y)^\beta)$ on setting $ x=\alpha^{-1}(\alpha y)^\beta$. 

Other integration techniques can be similarly used to create new distributions.
For example, Azzalini's method (\eg Azzalini and Capitanio, 2018) of transforming a symmetric pdf $f(x)$ to $2w(x)f(x)$,
where $w(-x)=1-w(x)$, can be thought of as a trick to simplify asymmetric integrands of type $w(x)f(x)$ used in reverse. 
These reflections prompt the thought that since a probability density function (pdf) must integrate to unity, all the `tricks' used to simplify unknown integrals,
if used in reverse, can be used to generate more complex integrands (pdfs). 

A method for simplifying integrands is integration by parts (IBP), and the creation of new distributions using this method was introduced by Baker (2019).
This article further explores the use of IBP and its discrete analogue, summation by parts (SBP) to generate new distributions. 

Before giving the general methodology, we show the power of the method with an example, starting with the exponential distribution.
This is a special case of a distribution given in Baker (2019).
Write the exponential pdf $f(t)=\alpha\exp(-\alpha t)$ for $\alpha > 0$ and  $T \ge 0$ as $f(t)=-uv^\prime$, 
where $u=\alpha t^{1/2}, v^\prime=-\exp(-\alpha t)/t^{1/2}$; this is just one of many choices of $u$ and $v$ that could be made.

Then $v(t)=\int_t^\infty x^{-1/2}\exp(-\alpha x)\dee x$, and 
on evaluating the integral by changing variable to $y$ where $\alpha x=y^2/2$ we see that $v(t)=\alpha^{-1/2} 2\sqrt{\pi}\Phi(-\sqrt{2\alpha t})$, where $\Phi$ is the normal distribution function.
We have that $u(0)=0, v(\infty)=0$, and hence using the method of parts, integrating $v^\prime$ and differentiating $u$, the integrand (the new pdf) is 
\[g(t)=\frac{\alpha \sqrt{\pi}\Phi(-\sqrt{2\alpha t})}{\sqrt{\alpha t}}.\]
This is a 1-parameter distribution, that can be used  for modelling lifetimes when the hazard function initially high, and the pdf and hazard function
are shown in figure \ref{fig:ca}. Further details are given in appendix A.

\begin{figure}
\centering
\makebox{\includegraphics{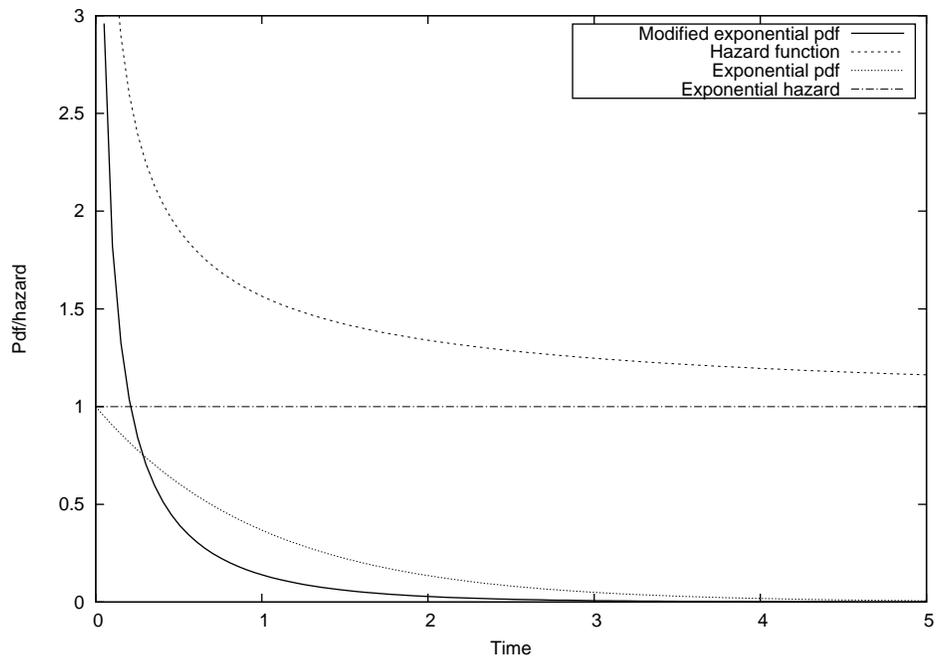}}
\caption{\label{fig:ca}The modified exponential pdf for $\alpha=1$ and its hazard function, along with the pdf and hazard function for the exponential distribution.}
\end{figure}

This example shows how potentially useful and fairly tractable new distributions can be generated quite easily using IBP, and illustrates the point that the pdf is often a special function.

The next section gives the general methodology and some general properties of the new distributions for the continuous case, and then some promising distribution classes are introduced.
For discrete distributions, the method of summation by parts can be used analogously and gives broadly similar results, discussed in section \ref{summ}.
\section{General theory\label{sec:genprop}}
\subsection{Deriving probability density functions}
Some notation is now introduced and the basic idea described more formally.
The method gives a transformation of the integrand (pdf) leading to a new pdf, and the two pdfs can conveniently be called `L' and `R', where the R form stochastically dominates the L form, so that the probability mass is further to the right.
The original distribution is sometimes referred to here as the `base' distribution.

Let the support of a distribution be $(x_l, x_h)$, which covers distributions defined on the whole real line, survival distributions, and doubly-bounded distributions.
Write the R-pdf as $f(x)=-u(x)v^\prime(x)$, where $u$ is a positive monotone increasing function so that $u^\prime \ge 0$, with $u(x_l)=0$, and $v$ a positive decreasing function so that $v^\prime(x) \le 0$, with $v(x_h)=0$.
Then applying integration by parts,
\begin{equation} \int_{x_l}^{x_h} f(x)\dee x=1=[-u(x)v(x)]^{x_h}_{x_l}+\int_{x_l}^{x_h} u^\prime(x)v(x) \dee x.\label{eq:parts}\end{equation}
The first term on the right vanishes, giving the L-pdf 
\begin{equation}g(x)=u^\prime(x)v(x).\label{eq:g}\end{equation}
When $u(x_l) > 0$ and $v(x_h) > 0$, (\ref{eq:parts}) yields a pdf 
\begin{equation}g(x)=u^\prime(x)v(x)/\{1+u(x_h)v(x_h)-u(x_l)v(x_l)\}.\label{eq:full}\end{equation}
Here by construction $v(x_h)=0$.
In Baker (2019) $u, v$ are crafted so that the simpler form
(\ref{eq:g}) could be used.
One can of course transform the integrand by going $\text{L}\rightarrow \text{R}$ from the base distribution 
instead of $\text{R}\rightarrow \text{L}$ as shown above,  to obtain $f(x)$ from $g(x)$ by integrating $u^\prime$ and differentiating $v$.

\subsection{Distribution functions and Stochastic dominance}
Let the distribution functions corresponding to $f(x), g(x)$ be $F(x), G(x)$ respectively.
Write the survival distributions $1-F(x)=\bar{F}(x)$ \etc
Then integrating by parts again
\begin{equation}G(x)=\int_{x_l}^x u^\prime(y)v(y)\dee y=[u(y)v(y)]^x_{x_l}-\int_{x_l}^x u(y)v^\prime(y)\dee y=u(x)v(x)+F(x),\label{eq:distfun}\end{equation}
and $\bar{G}(x)=\bar{F}(x)-u(x)v(x)$.
The L-distribution has probability mass shifted to the left, and the R-distribution dominates it stochastically, \ie $\bar{F}(x) > \bar{G}(x)$.
The mean can be written as $\text{E}_g(X)=\text{E}_f(X)-\int_{x_l}^{x_h}u(y)v(y)\dee y$.
The R-distribution may or may not dominate when (\ref{eq:full}) is used, and 
\begin{equation}G=\frac{F(x)+u(x)v(x)-u(x_l)v(x_l)}{1-u(x_l)v(x_l)}\label{eq:gplus}\end{equation}

The shifting of probability mass can be better understood by tagging the mass at $x_0$ using a Dirac delta-function.
Then we apply IBP to $-u(x)v^\prime(x)\delta(x-x_0)$, to obtain $-u^\prime(x)v^\prime(x_0)$ for $x \le x_0$, else zero.
Hence a unit probability mass at $x_0$ has been smeared out onto the range $(x_l,x_0)$ with pdf $u^\prime(x)/u(x_0)$, distribution function $u(x)/u(x_0)$.
Going from the L to the R distribution, the probability mass is smeared out over $(x_0,x_h)$ with survival function $v(x)/v(x_0)$.
\subsection{Random numbers}
The above suggests a method of generating random numbers from the L-distribution, given that we can generate them from the R-distribution $F$. Let $X$ be a random number from $F$
and $V$ a uniform random number. Then a random variable $Y$ from $G$ satisfies $u(Y)/u(X)=V$ so that $Y=u^{-1}\{u(X)V\}$. If $u$ can be inverted, this is a simple way to generate random numbers from the L-distribution.
\subsection{Families of distributions}
One usually seeks to generalize a pdf so that a family of distributions can be generated indexed by one or more parameters, with the original distribution as a special case.
The freedom available in choosing $u$ and $v$ makes this possible because 
if $v^\prime$ decreases infinitely fast, the integral $u(x)v(x)=-\int_x^{x_h}u(x)v^\prime(y)\dee y$ is zero, 
and $g(x)=f(x)$ for some value of the indexing parameter ($\lambda^{-1}=0$ in the examples given later).
This makes large-sample inference easier, as one can \eg use twice the log-likelihood increase on adding the new parameter in a chi-squared test of whether the 
model fit to a dataset has improved. It also enables parametric tests of stochastic dominance: appendix B gives a brief discussion.

The next section discusses particular classes of distributions generated using IBP.
\section{Some distribution classes\label{seclass}}
\subsection{Using a function of the distribution function for $u$}
With $u$ a function of $F$, $G$ will be a function of $F$, so this yields general families of distributions.
The most interesting classes have $G\rightarrow F$ as the parameter $\lambda\rightarrow\infty$.
In general, when $\dee u(x)/\dee F(x) \rightarrow \infty$ as $F \rightarrow 0$, the pdf $g(x)$ will be infinite at the lower limit $x_l$,  if $f(x_l)$ is finite.
\subsubsection{$u(x)=F^\lambda(x)$}

Moving $\text{R}\rightarrow\text{L}$, one choice is to take $u(x)=F^\lambda(x)$, where $\lambda > 0$, so that $v^\prime(x)=f(x)/F^\lambda(x)$.
Then $v(x)=\int_{F(x)}^1\dee y/y^\lambda=\frac{1-F^{1-\lambda}(x)}{1-\lambda}$.
Hence the L-distribution has pdf
\begin{equation}g(x)=\lambda\left(\frac{F^{\lambda-1}(x)-1}{1-\lambda}\right)f(x),\label{eq:f}\end{equation}
and distribution function 
\begin{equation}G(x)=\frac{F^\lambda(x)-\lambda F(x)}{1-\lambda}.\label{eq:glambda}\end{equation}
This is a negative mixture of the original pdf with weight $-\lambda/(1-\lambda)$ and the top $\lambda$-th order statistic with weight $1/(1-\lambda)$.
As $\lambda\rightarrow\infty$,  $G(x)\rightarrow F(x)$ and as $\lambda\rightarrow 0$, $G(x)\rightarrow 1$, \ie the probability mass is zero at above the lower limit.
When $\lambda=1$, from (\ref{eq:f}) or directly, $g(x)=-\ln(F)f(x)$ and $G(x)=\{1-\ln(F(x))\}F(x)$.

Expanding $F=1-\bar{F}$ in a Taylor series for small $\bar{F}$, we see in the right tail that $\bar{G}(x) \simeq (\lambda/2)\bar{F}^2(x)$, showing the pulling in of the right tail that follows from the movement of probability mass to the left.
The hazard function is $2f(x)/\bar{F}(x)$, twice that of the base distribution.

One may wonder how this is possible, given that $G \rightarrow F$ as $\lambda \rightarrow \infty$. The solution to the paradox is that the tail where the hazard
is double that of the base distribution occurs at larger and larger $x$ as $\lambda\rightarrow\infty$.
For small $F$, when $\lambda > 1$ we have that the hazard function $h(x) \simeq \lambda f(x)/(\lambda-1)$, so that the hazard and pdf is similar to the base distribution,
but both are larger by a factor $\lambda/(\lambda-1)$. For $\lambda < 1$, we have that $h(x) \simeq \lambda F(x)^{\lambda-1}f(x)/(1-\lambda)$, so that if $f(x)$ is finite,
the hazard would be infinite.

Going $L\rightarrow R$, set $v(x)=\bar{G}^\lambda(x)$, so $u^\prime=g\bar{G}^{-\lambda}$ and $u=\frac{1-\bar{G}^{\lambda-1}}{1-\lambda}$ and
$\bar{F}=\frac{\bar{G}^\lambda-\lambda \bar{G}}{1-\lambda}$. This is a negative mixture that is longer-tailed than was the L-distribution.
When $\lambda=1$, $\bar{F}=(1-\ln(\bar{G}))\bar{G}$.
\subsubsection{$u=\exp(\lambda F)$}

On taking $u(x)=\exp\lambda F(x)$ and using (\ref{eq:gplus}),
\begin{equation}G=\frac{\lambda F-\exp(-\lambda(1-F))+\exp(-\lambda)}{\lambda-1+\exp(-\lambda))}.\label{eq:gexp}\end{equation}
Since $H=\frac{\exp(-\lambda(1-F))-\exp(-\lambda)}{1-\exp(-\lambda)}$ is a distribution function, $G$ is the negative mixture 
\[G=\frac{\lambda}{\lambda-(1-\exp(-\lambda))}F-\frac{1-\exp(-\lambda)}{\lambda-(1-\exp(-\lambda)}H.\]
When $F$ is exponential, $H$ is a truncated extreme-value distribution.

AS $\lambda\rightarrow\infty$, (\ref{eq:gexp}) yields $G\rightarrow F$, and as $\lambda\rightarrow 0$, $\bar{G} \rightarrow \bar{F}^2$, \eg
a component lifetime distribution tends to the distribution of the lifetime of a series system where either of two components must fail for a failure of the system.

In the tail where $\bar{F}$ is small, expanding (\ref{eq:gexp}) shows that the hazard function is twice that of the base distribution (which is true for all $x$ as $\lambda\rightarrow 0$).

In the left tail, 
\[G \simeq \frac{(1-\exp(-\lambda))\lambda F}{\lambda-1+\exp(-\lambda)},\]
\ie the pdf is scaled up from the base distribution by a factor not exceeding 2.
For survival distributions as base distribution, this distribution therefore tends to give an increasing hazard function.

Other increasing functions of $F$ can also be used, \eg $u=\exp(\lambda F)-1$.
This gives 
\[G=F+(\exp(\lambda F)-1)\ln\left(\frac{1-\exp(-\lambda)}{1-\exp(-\lambda F)}\right)/\lambda,\]
which is messy.
Any increasing function can be used for $u$, leading to a vast number of possible (and messy) distributions.
\subsection{Using the transformation $u=t^\lambda$ for survival and doubly bounded distributions\label{sec:x}}
When $u$ is a power of the random variable (time $t$), more specific results follow. 
All lifetime distributions must have a scale factor $\alpha$. Without loss of generality, we set $\alpha=1$ for now.
The base pdf may include a power of $t$. Let this term be $t^{\beta-1}$, where $\beta$ is positive for Weibull and gamma distributions, and $\beta=1$ for the lognormal distribution.
Take $u=t^{\beta-1+\lambda}$. Then from (\ref{eq:g})
\begin{equation}g(t)=(\beta-1+\lambda)t^{\beta-2+\lambda}v(t).\label{eq:gee}\end{equation}
We must have that $\beta-1+\lambda > 0$ so that the pdf is positive, hence $\lambda > 1-\beta$. 

Baker (2019) discusses this case in detail and defines reliability growth as
\begin{equation}\xi=1/(\beta-1+\lambda).\label{eq:defxi}\end{equation}
The mean lifetime is $\text{E}_g(T)=\frac{\beta-1+\lambda}{\beta+\lambda}\text{E}_f(T)$, so the proportional increase in expected lifetime on eliminating inferior items is
\[\Delta\text{E}(T)=\frac{\text{E}_f(T)-\text{E}_g(T)}{\text{E}_g(T)}=\frac{1}{\beta-1+\lambda}=\xi,\]
and the proportional decrease in lifetime caused by inferior items is $1/(\beta+\lambda)=\xi/(\xi+1)$. 
As $\lambda$ increases and the base distribution is regained, $\Delta\text{E}(T)$ goes to zero.

\section{Special cases of continuous distributions}
This section briefly discusses specific distributions; 
a great variety of distributions can be generated, of which these are only a small sample.

With the exponential distribution $g(t)=\exp(-t)$ as base, shifting right using the $\bar{F}^\lambda(t)$ transformation from section \ref{seclass}, we have $\bar{G}=\exp(-x)$ and $\bar{F}(x)=\frac{\exp(-\lambda x)-\lambda\exp(-x)}{1-\lambda}$.
This is a 2-parameter phase-type distribution that occurs \eg as the prevalence of intermediate radioactive decay products. The hazard function increases linearly with slope $\lambda$ and 
levels off at $\min(1,\lambda)$. The moments are tractable, \eg the mean is $1+\lambda^{-1}$.

Using the $t^\lambda$ transformation to shift right yields $f(t)=\frac{ x\exp(-x)}{\lambda+1}+\frac{\lambda \exp(-x)}{\lambda+1}$, a mixture of exponential and gamma distributions.

Baker (2019) explored left-shifting distributions such as gamma and Weibull. These are special cases of the Stacy distribution.
The Stacy or generalized gamma distribution  has pdf 
\[f(t)=\frac{\alpha\gamma(\alpha t)^{\beta\gamma-1}\exp(-(\alpha t)^\gamma)}{\Gamma(\beta)},\]
where $\alpha > 0, \beta > 0, \gamma > 0$ and $\Gamma$
is the gamma function. It includes gamma, Weibull and lognormal distributions as special cases. Going $L\rightarrow R$ with $v(t)=\exp(-(\alpha t)^\gamma)/(\alpha t)^\lambda$
yields a simple mixture of Stacy distributions, but going $R\rightarrow L$ is more productive.
In integrating $v^\prime(t)$, the incomplete gamma function will be needed, defined as
\begin{equation}\Gamma(a;t)=\int_t^\infty x^{a-1}\exp(-x)\dee x,\label{eq:incgam}\end{equation}
where $a > 0$.
The gamma function itself is $\Gamma(a)=\Gamma(a;0)$.
If $a \le 0$, $\Gamma(a;t)$ is still defined but cannot be computed using software that computes the incomplete gamma function.

The resulting pdf is
\begin{equation}g(t)=\frac{\alpha(\lambda+\beta\gamma-1)(\alpha t)^{\beta\gamma+\lambda-2}\Gamma(\frac{(1-\lambda)}{\gamma};(\alpha t)^\gamma)}{\Gamma(\beta)},\label{eq:surv}\end{equation}
and the survival function is
\[\bar{G}(t)=\bar{F}(t)-\frac{(\alpha t)^{\beta\gamma+\lambda-1}\Gamma(\frac{(1-\lambda)}{\gamma};(\alpha t)^\gamma)}{\Gamma(\beta)}.\]
The relation between the moments for L and R-distributions for the $u=t^\lambda$ case was given in section \ref{sec:x}.

The beta distribution is the most common doubly-bounded distribution. Here the pdf $f(x)=B^{-1}(\alpha,\beta)x^{\alpha-1}(1-x)^{\beta-1}$, where $B$ denotes the beta function.
Taking $u=x^{\alpha+\lambda-1}$, we have that $v(x)=B(\alpha,\beta)^{-1}\int_x^1 y^{-\lambda}(1-y)^{\beta-1}\dee y+c$, where $c \ge 0$.
The constant of integration had to be zero for survival distributions as $u(\infty)=\infty$, but here does not.
Hence 
\[v(x)=B(1-\lambda,\beta;x)/B(\alpha,\beta)+c,\]
where $B(1-\lambda,\beta;x)$ is the complement of the unregularized incomplete beta function.
This yields the pdf 
\[g(x)=u^\prime v/(1+c)=\frac{(\alpha+\lambda-1)x^{\alpha+\lambda-2}v(x)}{1+c}.\]
The distribution function is
\[G(x)=\frac{F(x)+x^{\alpha+\lambda-1}(B(1-\lambda,\beta;x)/B(\alpha,\beta)+c)}{1+c},\]
where $F(x)=(1-B(\alpha,\beta;x))/B(\alpha,\beta)$.

On integrating $gx^n$ by parts, we have that
\[\text{E}_g(X^n)=\frac{\alpha-1+\lambda}{1+c}\{\frac{c}{\alpha-1+\lambda+n}+\frac{\text{E}_f(X^n)}{\alpha-1+\lambda+n}\}.\]

We have that $g(1)=(\alpha+\lambda-1)c/(1+c)$.
This 4-parameter distribution allows a much more flexible pdf than the beta distribution. The transformation can of course also be applied to $1-X$
by changing $X \leftrightarrow 1-X$ in the transformed distribution. Note that the beta distribution is label-invariant ($1-X$ also follows a beta distribution)
but the transformed distribution is not.

Random numbers can be generated as follows:
\begin{enumerate}
\item With $U$ a uniform r.v., if $U < c/(1+c)$ , generate $Z=V^{1/(\alpha-1+\lambda)}$, where $V$ is a uniform r.v.;
\item if $U \ge c/(1+c)$ generate $X$, a r.v. from the parent beta distribution ;
\item Then $Y=W^{1/\alpha+\lambda-1}X$, where $W$ is a uniform r.v.
\end{enumerate}
To keep generation of random numbers to a minimum, $V$ and $W$ could be generated by affine transformations on $U$.

On the whole real line, the normal distribution is the most important by far, and has been skewed in many ways, \eg Azzalini and Capitanio (2018).
Taking $u(x)=\exp(\lambda x)$ leads to the lagged-normal or normal-exponential distribution (\eg  Johnson {\em et al} 1995). This is exponential in the tail,
and can also be derived as the sum of normal and exponential random variables, so is not new. The $u=F^\lambda$ transformation yields the class (\ref{eq:glambda})
which with $F(x)=\Phi(x)$, where $\Phi$ is the normal distribution function, yields a skew distribution where $G(x)=(\Phi(x)^\lambda-\lambda \Phi(x))/(1-\lambda)$.
The moments cannot be found simply.
\section{Data fitting}
The method generates vast numbers of distributions, some of which are new.
The aim of the analysis was merely to show by an example that the new distributions can be useful.

In their book on survival analysis, Klein and Moeschberger (2003)  use a dataset of days to death of 863 kidney-transplant patients whose transplants were performed at the Ohio State University Transplant Center 
between 1982 and 1992. Available covariates are age, gender and white/black. Only 140 patients' survival times were not censored.
We discretized age into 5 bands: 1-16, 17-32, 33-48, 49-64 and 65+ and fitted the distribution from (\ref{eq:surv}) with $\beta=1$, \ie a modified Weibull distribution.
An accelerated time model was used (\eg Chiou {\em et al}, 2014), so that $\alpha=\alpha_0\exp({\boldsymbol \eta}^T{\bf X})$, where ${\bf X}$ is a vector of covariates and ${\boldsymbol \eta}$
a vector of regression coefficients. The model parameters are then the $\boldsymbol \eta$, $\alpha_0$ the baseline time-scale, $\gamma$ the Weibull shape parameter
and $\lambda$, reparameterized so that $\xi$ as defined in (\ref{eq:defxi}) is used instead of $\lambda$. The fit was done by maximum-likelihood in a purpose-written fortran program.
The incomplete gamma function for negative argument is not available as a standard special function, and so was evaluated as an integral for $a \le 0$ in (\ref{eq:incgam}).
The NAG library routine D01AMF was used. This routine transforms the integration range to be finite and then integrates adaptively.
Similar routines exist on many platforms.

\begin{figure}
\centering
\makebox{\includegraphics{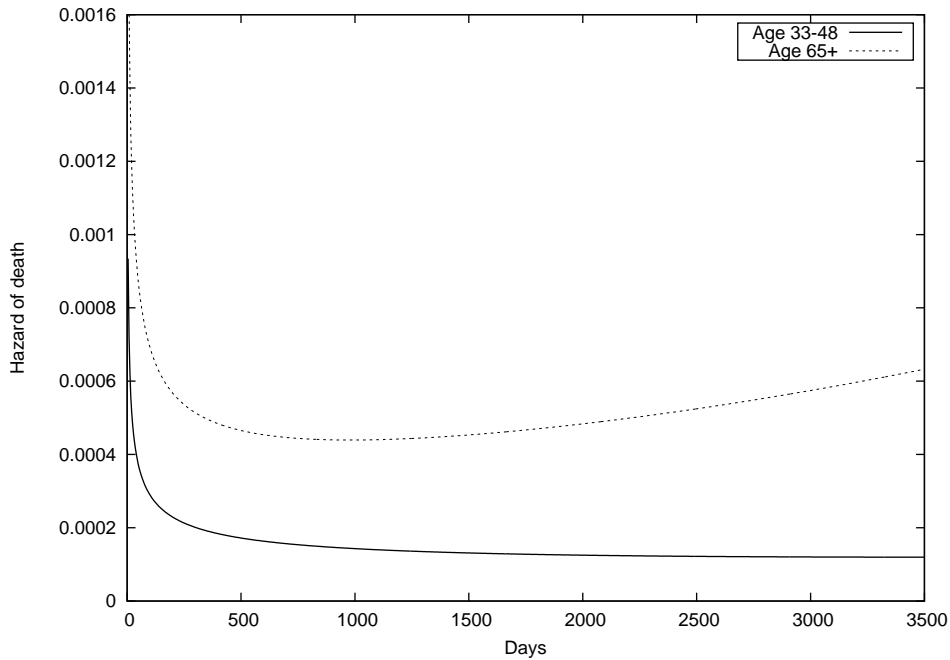}}
\caption{\label{fig:cb}Hazard function for the modified Weibull distribution fitted to the kidney transplant data of Klein and Moeschberger (2003)
for the 33-48 and 65+ age groups.}
\end{figure}

The fit to data cannot be shown because of the censoring, but
figure \ref{fig:cb} shows the hazard function of the fitted distribution for the central and oldest age-bands. From this it can be seen that
the hazard of death decreases steadily after transplant, effectively to a constant for the central age-band, but starts increasing again for the top age band.
This result is consistent with hazard plots using all patients produced by various methods by Klein and Moeschberger (2003).

The table shows fitted parameters and standard errors. 
Standard errors are large, sometimes very large, because little information comes from the censored survival times.
The baseline age group was 33-48 years.
It can be seen that the hazard of death increases rapidly from the lowest age band,
and then increases more slowly with age. Women have a slightly higher hazard of death than men, and people of colour higher than white, but these effects are not statistically significant.
Because $\gamma > 1$, the hazard of death would rise eventually for all age groups.
The value of 1.6 for $\xi$ means that, for a given set of covariate values, the expected lifetime could be increased by 160\% of its current value,
if all patients could be made to respond as well as the best ones. This suggests that much could potentially be gained by studying and improving survival rates.

\begin{table}[h]
\begin{tabular}{|l|c|c|} \hline
Parameter & Estimate& Standard error  \\ \hline
$\alpha_0$ &0.576$\times 10^{-4}$ &  0.261$\times 10^{-4}$   \\ \hline
Shape  $\gamma$ &3.315   &     4.683    \\ \hline
Mixing parameter $\xi$ &1.603 &      0.123    \\ \hline
Gender & 0.048 & 0.240   \\ \hline
Black/white &0.243  &    0.295\\ \hline
1-16 age band &-29.685     &  168.886\\ \hline
17-32 age band &-2.354     & 0.530     \\ \hline
49-64 age band &-.558      &0.266\\ \hline
65+ age band & 1.338  &    0.377 \\ \hline
\end{tabular}
\caption{Estimated parameters and standard errors for the fit to the kidney transplant dataset.}
\end{table}
\section{Discrete distributions from summation by parts\label{summ}}
A way to generalize discrete (count data) distributions,
is to use partial sum distributions.
Wimmer and Ma\v{c}utek (2012) give a general definition of partial-sum distributions, where 
for parent (base) distributions with support on the non-negative integers, with probabilities $p_0, p_1 \ldots$
the descendant probabilities $q_i$ are given by $q_i=\sum_{j \ge i}f(i,j)p_j$, where $f$ is a real function.
The simplest such distribution has $q_i=\mu^{-1}\sum_{j=i+1}^\infty p_j$ where  $\mu$ is the mean of the parent distribution; this is a discrete form of the length-biased distributions found in renewal theory. 
Wimmer and Altmann (2001) describe more general distributions of this type, and Johnson {\em et al} (2005) give a summary.

Here a methodology is given that yields classes of partial-sum distributions.
The two classes that have been explored in detail reduce to the parent distribution when a parameter $\lambda$ or $r \rightarrow \infty$.
We shall call them $r$-Poisson, $\lambda$-binomial, \etc
The r-class interpolates between a distribution and a variant of its length-biased form. 

The population mean is needed for statistical inference.
For example, one often wishes to regress the mean on some covariates. This is computationally more difficult (but still possible) when the mean cannot be readily computed in terms of model parameters.
The class of r-distributions has the advantage of tractable moments and so allows more straightforward inference.

The general methodology is next given, followed by the properties of the 2 classes of distributions
and example of a fit to data is given.

\subsection{The general pmf}
The mathematical technique of summation by parts (SBP) can be used to transform distributions into new ones. The SBP identity can be written
\begin{equation}-\sum_{i=m}^n u_i(v_{i+1}-v_i)=u_mv_m-u_nv_{n+1}+\sum_{i=m+1}^n v_i(u_i-u_{i-1})\label{eq:summ}\end{equation}
for any $u_i, v_i$.
The proof follows from the telescoping property of the sum $\sum_{i=m+1}^n (u_iv_{i+1}-u_{i-1}v_i)=u_nv_{n+1}-u_mv_{m+1}$.
Distributions are commonly defined on the integers $0$ to $n$, where  usually, $n=\infty$, the main exception being the binomial distribution. 
Let a `parent' probability mass function (pmf) $p_i=-u_i(v_{i+1}-v_i) \ge 0$, $u_i \ge 0$ and let $u_{i+1}-u_i \ge 0$ with $v_{n+1}=0$.
Since $p_i \ge 0$, $v_i \ge v_{i+1} \ge 0$.
Then from (\ref{eq:summ}) 
\begin{equation}q_i=\frac{ v_i(u_i-u_{i-1})}{1-u_mv_m}\label{eq:qdef}\end{equation}
is also a pmf.

Set $m=-1$ and set the pmf $p_{-1}=0$. 
Then from (\ref{eq:qdef}) for $i \ge 0$
\[q_i= \frac{v_i(u_i-u_{i-1})}{1-u_{-1}v_{-1}}.\]

Since $v_{i+1}-v_i=-p_i/u_i$ and $v_{n+1}=0$, it follows that $v_i=\sum_{j=i}^n p_j/u_j$. Since $p_{-1}=0$, $v_{-1}=v_0$.
The function $u_i$ can be chosen to be zero at $m=-1$, \eg $u_i=(i+1)^\lambda$ where $\lambda > 0$, or $r^i-1/r$ where $r > 1$.
In this case simply $q_i=(\sum_{j=i}^n p_j/u_j)(u_i-u_{i-1})$.
Otherwise
\begin{equation}q_i= \frac{(\sum_{j=i}^n p_j/u_j)(u_i-u_{i-1})}{1-u_{-1}v_0}.\label{eq:qf}\end{equation}
Using (\ref{eq:summ}) in this way, both parent and descendant distributions have the same upper limit $n$ and hence the same support.

Note that for $q_i \ge 0$ we require $u_{-1}v_0 < 1$. Since $v_0=\sum_{j=0}^n p_j/u_j < 1/u_0$ and $u_{-1} < u_0$, we have that $u_{-1}v_0 < 1$ as required.
If $u$ is indexed by a parameter $r$ so that $u_{i+1}/u_i \rightarrow \infty$ as $r\rightarrow\infty$, then from (\ref{eq:qf}) $q_i \rightarrow p_i$ as $r \rightarrow\infty$.
This is so because $v_i \rightarrow p_i/u_i$, $u_i-u_{i-1} \rightarrow u_i$, and $u_{-1}v_0 \rightarrow p_0(u_{-1}/u_0)\rightarrow 0$.
Hence the class of descendant distributions generalizes the  parent distribution. 
\subsection{The discrete distribution function and moments}
Denoting distribution functions for parent and descendant distributions by $F_k, G_k$ respectively, we have from (\ref{eq:summ}) that 
\begin{equation}G_k=\frac{F_k+u_kv_{k+1}-u_{-1}v_0}{1-u_{-1}v_0}.\label{eq:df}\end{equation}
Hence the distribution function is readily calculable once the $v_i$ are calculated and the pmf is known.
If $u_{-1}=0$, $G_k > F_k$ and the parent stochastically dominates the descendant distribution.

From (\ref{eq:qf}) the mean $\mu_g$ is
\begin{equation}\mu_g=\frac{\sum_{i=0}^n i(u_i-u_{i-1})\sum_{j=i}^n p_j/u_j}{1-u_{-1}v_0}.\label{eq:mean}\end{equation}
On reversing the order of summation,
\begin{equation}\mu_g=\frac{\mu-\sum_{i=1}^n(\sum_{j=0}^{i-1}u_j)p_i/u_i}{1-u_{-1}v_0}.\label{eq:mugsumm}\end{equation}
For the $r$-th non-central moment,
\begin{equation}\mu_g^{(r)}=\frac{\mu^{(r)}-\sum_{i=1}^n(\sum_{j=0}^{i-1}\{(j+1)^r-j^r\}u_j)p_i/u_i}{1-u_{-1}v_0}.\label{eq:moms}\end{equation}
\subsection{Random numbers}
A unit probability mass at $k$ for the parent distribution becomes a pmf $r_i=(u_i-u_{i-1})/(u_k-u_{-1})$ for $i \le k$, so random numbers can be generated by generating a random number $K$ from the parent distribution,then
choosing a random number $I$ with probability $r_i$. If $K$ is not too large, this can be done by generating a uniform random number $U$,
and choosing the smallest $i$ such that $\sum_{j=0}^i r_j > U$. This is the table lookup method (\eg Shmerling, 2013).

From (\ref{eq:df}) on setting $k=0$, we have that $q_0=(p_0+u_0v_1-u_{-1}v_0)/(1-u_{-1}v_0)$, so $q_0 > p_0$ when $u_{-1}=0$, which makes these distributions useful when the parent distribution (\eg Poisson) underfits the probability of no events occurring.

\subsection{The $u_i=(i+1)^\lambda$ distribution}
With the choice $u_i=(i+1)^\lambda$, $u_{-1}=0$ and the new system of distributions
\begin{equation}q_i=\{(i+1)^\lambda-i^\lambda\}\sum_{j=i}^n p_j/(j+1)^\lambda\label{eq:q2}\end{equation}
is obtained. 

As $\lambda\rightarrow\infty$, $q_i \rightarrow p_i$, and as $\lambda\rightarrow 0$, $q_0 \rightarrow 1$ and all the probability mass
is at zero.

The case where $\lambda=1$ is briefly mentioned in Johnson {\em et al} (2005).
In general, the moments are intractable, unless $\lambda=1$. There (\ref{eq:mugsumm}) enables the mean to be found,
and higher moments are also tractable.
Letting $p_j\rightarrow p_{j+1}$ and setting $p_0=0$ gives 
\[q_i=\{(i+1)^\lambda-i^\lambda\}\sum_{j=i+1}^\infty p_j/j^\lambda.\] 
This is a class of distributions that generalizes the Bissinger system of distributions (Johnson {\em et al} 2005, p. 509),
in which $\lambda=1$.

The choice $u_i=r^i$ where $r > 1$ is more tractable, and is discussed next.
\subsection{The $u_i=r^i$ distribution}
\begin{figure}
\centering
\makebox{\includegraphics{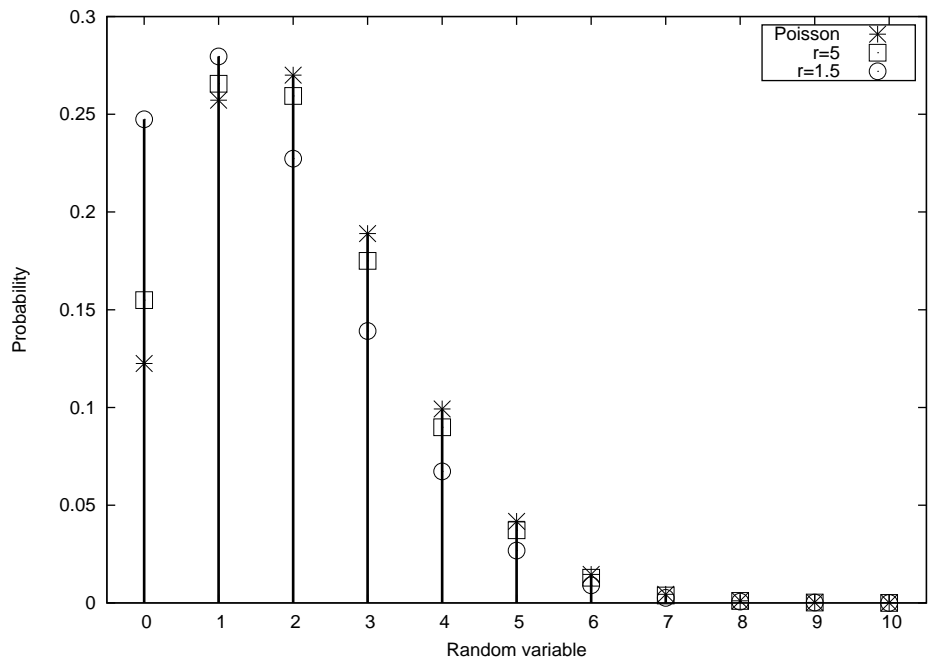}}
\caption{\label{dfiga}Probabilities for the Poisson distribution with $\mu=2.1$, and for the r-distribution with $r=1.5$ and $r=5$.}
\end{figure}

Here $u_i=r^i$ for $r > 1$. Define $H_i(x)=\sum_{j=i}^n p_jx^j$, then $v_i=H_i(1/r)$ and 
\begin{equation}q_i=\frac{H_i(1/r)r^i(1-1/r)}{1-H_0(1/r)/r},\label{eq:qh}\end{equation}
as $v_{-1}=1/r$.
Also, $G_i=\frac{F_i+r^iH_{i+1}(1/r)-H_0(1/r)/r}{1-H_0(1/r)/r}$.
The pmf for a Poisson parent distribution is shown in figure \ref{dfiga}.

The function $H_0$ is the probability generating function (pgf) for the parent distribution.
From (\ref{eq:qh}) on reversing the order of summation, the pgf $M(s)$ for the derived distribution is 
\[M(s)=\frac{r-1}{rs-1}\frac{\{sH_0(s)-H_0(1/r)/r\}}{1-H_0(1/r)/r}.\]
From this, the moments can be read off, \eg the mean $\mu_g$ is
\begin{equation}\mu_g=\frac{\mu+1}{1-H_0(1/r)/r}-\frac{r}{r-1}.\label{eq:mug}\end{equation}
The pgf $H_0$ is known in analytic form for the major discrete distributions, \eg for the Poisson, $H_0(1/r)=\exp(-(1-1/r)\mu)$.
Hence the moments of the descendant distribution can be found as functions of the parameters of the parent distribution.

Random numbers can be computed using the general method described earlier. This now particularizes to the following: generate a random number $K$ from the parent distribution.
Then compute 
\[m=\frac{\ln\{(r^K-r^{-1})U+r^{-1}\}}{\ln(r)}.\]
where $U$ is a uniform $[0,1]$ random number. Then the integer part of  $m+1$ is $M$, a r.v. from the descendant distribution.
This is the inverse-transform method (Shmerling, 2013).

As $r \rightarrow\infty$ we have that $q_i\rightarrow p_i$, while when $r\rightarrow 1$ the numerator and denominator of (\ref{eq:qh}) go to zero.
Applying L'Hospital's rule and using $\dee H_0(x)/\dee x|_{x=1}=\mu$, $q_i \rightarrow (p_i+\sum_{j=i+1}^n p_j)/(1+\mu)$, a partial sum distribution corresponding to the parent distribution.
\subsection{Long-tailed discrete distributions}
Going from $q$ to $p$ a new class of long-tailed distributions can be obtained. 
For doubly-bounded discrete distributions, one can use (\ref{eq:qf}), with $p_i \rightarrow p_{n-i}$.
In general, (\ref{eq:summ}) is used in reverse. 
For example, for the $\lambda$-distributions, set $v_i=(i+1)^{-\lambda}-(n+2)^{-\lambda}$ so $u_i-u_{i-1}=q_i/\{(i+1)^{-\lambda}-(n+2)^{-\lambda}\}$.
This yields $u_i=\sum_{j=0}^i q_j/\{(j+1)^{-\lambda}-(n+2)^{-\lambda}\}$ so that the pmf $p_i$ is given by
\[p_i=\{(i+1)^{-\lambda}-(i+2)^{-\lambda}\}\sum_{j=0}^i q_j/\{(j+1)^{-\lambda}-(n+2)^{-\lambda}\}.\]
As $n \rightarrow\infty$ this simplifies to
\[p_i=\{(i+1)^{-\lambda}-(i+2)^{-\lambda}\}\sum_{j=0}^i q_j(j+1)^{\lambda}.\]

For such distributions $p_0=(1-2^{-\lambda})q_0$, so $p_0 < q_0$, and the probability of zero events is reduced. The probabilities $p_i$ are easy to compute,
as an infinite sum is not required.

The $r$-distribution with $n=\infty$ is trivial: the pgf is simply
\[M(s)=\frac{r-1}{r-s}H_0(s),\]
the product of the parent pgf and the pgf of a geometric distribution with $\text{Prob}(X=k)=r^{-k}/(1-1/r)$. Hence the effect of SBP has been to
add a geometric random variable to the original.
\subsection{Other discrete distributions}
Each functional form for $u_i$ gives a new family of distributions, and if $u_i$ can increase arbitrarily fast with $i$,
the parent distribution can be regained, so the descendant distribution will generalize the parent. For example, with $u_i=r^i-1/r$ one obtains
\[q_i=r^i(1-1/r)\sum_{j=i}^n p_j/(r^j-1/r),\] as $u_1=0$. As $r \rightarrow 1$, $q_i\rightarrow \sum_{j=i}^n p_j/(j+1)$. 
In general, the mean cannot be simply expressed.
Clearly, more parameters than one can be introduced, and
a host of new distributions created. In general, they will be quite messy, for example when $u_i=r^{\sqrt{i+1}}$. However, there may well be other functional forms for $u_i$
that give attractive distributions.

An interesting curiosity is a distribution obtained by setting $u_i=\prod_{s=1}^t(i+s)$, where $t$ is a positive integer.
This has $u_{-1}=0$, and $u_i-u_{i-1}=t\prod_{s=1}^{t-1}(i+s)$. The pmf is therefore
\[q_i=t\prod_{s=1}^{t-1}(i+s)\sum_{j=i}^n p_j/\prod_{k=1}^t(j+k).\]
The mean is tractable, and from (\ref{eq:mean}), $\mu_g=t\mu/(t+1)$.
From (\ref{eq:moms}) the variance is
\[\sigma^2_g=\frac{t\sigma^2}{t+2}+\frac{t\mu^2}{(t+1)^2(t+2)}+\frac{t\mu}{(t+1)(t+2)},\]
where $\sigma^2$ is the variance of the parent distribution.
As $t \rightarrow\infty$ the parent distribution is regained. This distribution generalizes the $t=1$ case mentioned briefly in Johnson {\em et al} (2005).
\section{Inference}
We usually wish to study the effect of covariates on the mean $\mu_g$.
The mean $\mu_g$ is therefore predicted as $\mu_g=\mu_0\exp({\boldsymbol \beta}^T{\bf X})$, where ${\bf X}$ is a vector of covariates. 
For the example, for the $r$-distributions, the model parameter $\mu$ was found for each case from (\ref{eq:mug}), using Newton-Raphson iteration, which usually converged in 4 iterations.
The pgf $H_0(x)=\{1+\mu\alpha(1-x)\}^{-1/\alpha}$ for the negative binomial distribution. 
The probabilities $q_i$ were then computed and fits were by maximum likelihood.

Note that for the $\lambda$-distributions the mean is not readily calculable in terms of parent distribution parameters.
Model fitting is still possible however, because the mean $\mu_g$ can be computed using the probabilities $q_i$ taken up to some large cutoff value of $i$,
and the Newton-Raphson iteration still done for $\mu$.
\section{Example}

Hilbe (2011) fits negative binomial distributions to a number of datasets. One is the `affairs' dataset, with 601 observations
from Fair (1978), reporting counts of extramarital affairs over a year in the USA. Table \ref{tab1} shows results of fitting the Poisson and NB distributions, and on adding
the $r$ parameter, and table \ref{tab2} shows the covariates and the regression coefficients.
The mean quoted is for the average values of the covariates.
\begin{table}[h]
\begin{tabular}{|l|c|c|} \hline
Model & -$\ell$ & AIC\\ \hline
Poisson & 1426.8 & 2871.54 \\ \hline
Poisson+$r$& 1126.62&2273.24 \\ \hline
Negative Binomial & 728.10& 1476.20 \\ \hline
Negative Binomial + r& 711.45 & 1444.91 \\ \hline
\end{tabular}
\caption{Model fits for models of the extramarital affairs data of Fair (1978), showing minus the log-likelihood and the Akaike Information Criterion\label{tab1}.}
\end{table}

\begin{table}[h]
\begin{tabular}{|l|l|l|} \hline
Variable & Coeff & p\\ \hline
Mean & 1.873 (.388) & - \\ \hline
$\alpha$ & 3.02 (.155) & - \\ \hline
$r$ & 1 & - \\ \hline
Gender & -.064 (.265) & .81 \\ \hline
Age & -.0229 (.0188) & .226 \\ \hline
Years married & .107 (.0355) & .0024 \\ \hline
Children?& .113 (.307) & .366 \\ \hline
Religious?& -.415 (.100) & .000036 \\ \hline
Educ. level & -.000610 (.0560) & .991 \\ \hline
Occupation & .0737 (.0801) &.920 \\ \hline
Rating & -.447 (.099) & .0000003 \\ \hline
\end{tabular}
\caption{Fitted parameter values for the NB+r model of the extramarital affairs data.
Standard errors are given in parentheses. The last column is the p-value for a test that the regression coefficient is zero.
\label{tab2}}
\end{table}
Clearly, adding the $r$ parameter to the negative binomial model has improved the fit significantly. The model has gone to the limiting case of the partial-sum distribution where $r \rightarrow 1$.
The conclusions about the effect of covariates agree broadly with Hilbe's analysis.
The significant predictors are self-rating of the marriage from unhappy to happy on a scale of 1-5, degree of religiosity on a scale from 1 to 5 (anti to very)
and years of marriage. Religious people with a happy marriage who have not been married long have fewest affairs.

This example shows how the new distributions can improve model fit and so enable better inference.
Computations were done with a purpose-written fortran program, that used the Numerical Algorithms Group (NAG) function minimisers.

\section{Conclusions}
Integration by parts is a general method that yields a cornucopia of distributions, only a few of which have been explored here.
Some are new, while some have been derived by other methods and now have a new characterization. This could be useful \eg in generating random numbers.

IBP is an addition to other general methods for modifying distributions, such as transforming variables. Summation by parts can be used similarly for discrete distributions,
and its usefulness may be relatively greater, as the variety of methods for generating new distributions from old is more limited in the discrete case.

One can shift the probability mass of a distribution left or right. Shifting it left is useful when dealing with failure-time distributions in reliability,
where the left-shifted distributions can reproduce the bathtub-shaped hazard functions sometimes seen in practice (\eg Lai and Xie, 2006). For discrete distributions, 
the left shifted distributions have a higher probability of zero events occurring. This zero-inflation can substantially improve the model fit to data.
Shifting probability mass right gives long-tailed distributions, which are needed in many areas, \eg finance.

Comparing modifying distributions by integration by parts to using the transformation of variables method, one could sum up as follows:
\begin{enumerate}
\item transformation gives a simple pdf, whereas IBP gives the pdf as an integral, which may be simple or may be a special function;
\item the distribution function for the transformation of variables method may not be tractable, but for IBP it is easily derived once the pdf is known
and is simply related to it by (\ref{eq:distfun});
\item the IBP method has the further property of yielding a new distribution that either stochastically dominates the original, or is dominated by it.
This property can be used in a test of stochastic dominance as described in appendix B.
\item The IBP method changes the skewness of distributions, but cannot change tail length whilst leaving skewness unchanged, as does \eg the arcsinh transformation
applied to the normal distribution.
However, IBP could be applied twice \eg $\text{R}\rightarrow \text{L}$ and then $\text{L}\rightarrow \text{R}$ to do this by lengthening each tail in turn.
\end{enumerate}

The methodology does not apply to circular distributions, but it is possible that an analogous method could be developed here.
Possible multivariate applications have not been mentioned. Clearly, with a bivariate distribution one can shift the probability mass of $X$ given $Y$,
to obtain a new distribution with a different marginal distribution for $X$, but the original marginal distribution for $Y$.
This type of procedure will not lead to a copula, because the marginal distribution of $X$ has changed, but such distributions could still be useful.

It is possible that other integration techniques, such as contour integration, could also be profitably used.
Much further work could be done in modifying distributions using partial integration, either developing the general methodology, or using it to derive more new distributions.
Fast computation of the resulting special functions such as (\ref{eq:incgam}) is also needed.

The mathematical method of summation by parts enables the creation of new families of discrete distributions that generalize existing distributions.
Any increasing function defined on non-negative integers (strictly, also on -1) gives a family of distributions, and two such functions, $u_i=(i+1)^\lambda$ and $u_i=r^i$ have been discussed.
The latter gives a relatively tractable distribution, with moments that can be found analytically in terms of the parameters of the parent distribution
when the pgf can be found analytically.

The distributions can be used for modelling data and statistical inference, and could give a useful sensitivity analysis when a model such as negative binomial has been fitted.
The class of long-tailed distributions of the form $q_i=\sum_{j \le i}f(i,j)p_j$ is completely new.

Further work could include derivation of new classes of distribution using other functions for $u_i$, and more detailed exploration of their properties.
Bivariate distributions have not been considered here, but new bivariate distributions could be derived from old by applying SBP to $X$ for each level of $Y$.

\clearpage
\section*{Appendix A: the modified exponential distribution}
The distribution is a special case of the modified Stacy distribution where $\beta=\gamma=1$, $\lambda=1/2$.
Some properties of this distribution are now given without proof, for completeness, and to show the tractability of the distribution.
The pdf is initially infinite, and is exponential in the tail.
The survival function is
\[\bar{G}(t)=\exp(-\alpha t)-2\sqrt{\pi}\Phi(-\sqrt{2\alpha t})(\alpha t)^{1/2}.\]
The hazard function $h(t)=g(t)/\bar{G}(t)$ decreases from infinity at $t=0$ to a constant $\alpha$, at large $t$. 

The moments are $\text{E}(T^n)=\alpha^{-n}n!/(1+2n)$, giving $\text{E}(T)=1/3\alpha$, $\text{E}(T^2)=(2/5)\alpha^{-2}$,
$\text{var}(T)=(13/45)\alpha^{-2}$. The coefficient of variation is thus $\sqrt{13/5} \simeq 1.61245$, compared with 1 for the exponential distribution.
Random numbers are generated as $T=-U^2\ln(V)/\alpha $, where $U, V$ are uniformly distributed on $[0,1]$.

The third moment about the mean is $502/945$, giving skewness of $\frac{502\sqrt{5/13}}{91} \simeq 3.42118$,
larger than 2 for the exponential distribution. The excess kurtosis is $17 \frac{145}{1183} \simeq 17.12257$, compared to 6 for the exponential distribution.

\section*{Appendix B: Stochastic dominance tests}
Given two samples  of returns from investments, one might wish to test for (first order) stochastic dominance of investment B by investment A (\eg Barrett and Donald, 2003).
A possible test is to fit a parametric form to the sample, with the indexing parameter fixed at zero (giving the R-distribution) for sample A,
and floating for sample B. The chi-squared test or other test that $\lambda^{-1} \ne 0$ can then be used to test whether investment A dominates B stochastically.
An exact test would compute $\hat\lambda$, then permute observations between the two groups A and B, refitting to generate the null distribution of $\hat\lambda$.
The p-value is the proportion of permuted samples for which $\hat\lambda$ is no greater than the observed value.
This test can also be done with  discrete distributions of gains from investment/gambling.
Note that unlike the description here, $H_0$ is often taken that there is stochastic dominance. This area of inference is complex and not the main subject of this article, but clearly good parametric tests could be devised.

\end{document}